\documentclass[11pt, reqno]{amsart}
\usepackage{diagbox, amsmath, amsthm, amscd, amsfonts, amssymb, graphicx, color, pgf, tikz}
\usepackage{booktabs}
\usepackage[bookmarksnumbered, colorlinks, plainpages]{hyperref}
\usepackage[all,cmtip]{xy} 
\usepackage{mathrsfs}
\usetikzlibrary{arrows}
\usepackage{amsmath}
\usepackage[autostyle]{csquotes}
\makeatletter \oddsidemargin.9375in \evensidemargin \oddsidemargin
\newtheorem{theorem}{Theorem}[section]
\newtheorem{lemma}[theorem]{Lemma}
\newtheorem{proposition}[theorem]{Proposition}
\newtheorem{corollary}[theorem]{Corollary}
\newtheorem{question}[theorem]{Question}
\theoremstyle{definition}
\newtheorem{definition}[theorem]{Definition}

\newtheorem{remark}[theorem]{Remark}
\numberwithin{equation}{section}

\begin{document}
\title[Finiteness and Ping-Pong Properties ]{Paradoxical decomposition of Group Actions and Configuration}

\author[M. Meisami, A. Rejali, M. Soleimani and A.Yousofzadeh]{M. Meisami, A. Rejali, M. Soleimani Malekan and A.Yousofzadeh}

\address{ Mahdi Meisami, Department of Mathematics, Ph.D student, University Of Isfahan , Isfahan 81746-73441,
Iran.}\email{m.meisami@sci.ui.ac.ir}
\address{Ali Rejali, Department of Mathematics, University Of Isfahan, Isfahan 81746-73441, Iran}\email{rejali@sci.ui.ac.ir}
\address{Meisam Soleimani Malekan, Department of Mathematics, University Of Isfahan , Isfahan 81746-73441,
Iran}\email{m.soleimani@sci.ui.ac.ir}
\address{Akram Yousofzadeh, Department of Mathematics, Mobarakeh Branch Islamic Azad University, Isfahan, Iran}\email{ayousofzade@yahoo.com}
\begin{abstract}
In this paper we study the notion of configuration for group actions. It is proved that some properties concerning configuration of groups can be extended for the case of group actions. The relationship between configuration and different types of well-known ping-pong lemma is also investigated.\end{abstract}

\maketitle

%%%%%%%%%%%%%%%%%%%%%%%%%%%%%%%%%%%%%%%%%%%%%%%%%%%%%%%%%%%%%%%%%%%%%%%%%
% Macros
%%%%%%%%%%%%%%%%%%%%%%%%%%%%%%%%%%%%%%%%%%%%%%%%%%%%%%%%%%%%%%%%%%%%%%%%%

\newcommand\sfrac[2]{{#1/#2}}

\newcommand\cont{\operatorname{cont}}
\newcommand\diff{\operatorname{diff}}

%%%%%%%%%%%%%%%%%%%%%%%%%%%%%%%%%%%%%%%%%%%%%%%%%%%%%%%%%%%%%%%%%%%%%%%%%
\section{\bf Introduction}
In 1929 von Neumann initiated the concept of amenable groups. He proved that no amenable group contains a free non-abelian subgroup \cite{von}.  The problem of holding the converse is known as the von Neumann conjecture \cite{day}. This conjecture disproved by Grigorchuk, presenting different sets of examples. Olshanskii constructed the first counterexamples \cite{olshanski}.  Counterexamples to the Burnside conjecture given by Adyan  are non-amenable torsion groups, which clearly contain no non-abelian free subgroups \cite{ady}.  Golod-Shafarevich groups are the first examples of torsion non-amenable residually finite groups. Other explicitly constructed examples can be found in \cite{sge} and \cite{monod}. However, there are families of groups with the property that each group in these families is either amenable or contains a free subgroup (see \cite{Tits} and \cite{monod} for example).

The action of a group $G$ on a non-empty set $X$ is a function $\cdot :G\times X\longrightarrow X$ such that the following conditions hold
\begin{itemize}
\item[(1)] $e\cdot x=x$ for each $x\in X$.
\item[(2)]
$g_1\cdot(g_2\cdot x)=(g_1g_2)\cdot x$ for all $g_1,g_2\in G, x\in X$.
\end{itemize}
 We say $G\curvearrowright X$ admits a paradoxical decomposition if there exist disjoint subsets $A_1,\dots,A_n,B_1,\dots,B_m$ of $X$ and elements $g_1,\dots ,g_n,h_1,\dots,h_m$ of $G$ such that  
$$X=\bigcup_{i=1}^n g_iA_i=\bigcup_{j=1}^m h_jB_j.$$
The Tarski number of $G\curvearrowright X$, which is denoted by $\tau(G\curvearrowright X)$ is the smallest number of pieces needed for a paradoxical decomposition of $G\curvearrowright X$. If this action admits no paradoxical decomposition, it is usual to write $\tau(G\curvearrowright X)=\infty$. If $G$ acts on itself by right multiplication, then we use the notation $\tau(G)$ for the Tarski number of $G\curvearrowright G$. It is easily checked that the smallest possible Tarski number of any group action is 4. By \cite{wagon,Theorem 4.5 & 4.8} Let $G\curvearrowright X$ be a group action the Tarski number of the action is four if and only if $G$ contains a non-abelian free subgroup $F$ such that the stablizers of points from $X$ in $F$ are cyclic. 

The group action $G\curvearrowright X$ is called amenable if there exists a positive finitely additive measure $\mu$ from powerset of $X$ to the $\rightarrow[0,\infty]$ such that 
\begin{itemize}
\item[(1)] $\mu(g\cdot E)=\mu(E)$ for each $E\subset X$ and $g\in G$.
\item[(2)] $\mu(X)=1$.
\end{itemize}

All groups and spaces are assumed to be discreye. The Tarski alternative states that a group $G$ is amenable if and only if it admits no paradoxical decomposition  \cite{wagon}. An analogous statement is valid for the group action $G\curvearrowright X$ (see \cite{wagon} or \cite{bhv} for more details). For more study of the configuration and amenability we refer the reader to the \cite{pier}, \cite{pat} and \cite{sur}.

The concept of configuration of groups can be extended for group actions. We shall study and investigate the properties of groups  and group actions related to this concept in section 2.

Let $G$ be a group acting on a set $X$. The well-known ping-pong lemma attributed to Felix Klein provides a criterion for determining when $G$ contains a free subgroup \cite[II. B]{harp1}. In section 3 we use configuration of groups and group actions and present some ping-pong like lemmas to characterize the certain properties of groups. In that section we also investigate the relationship between  Tarski number of group actions and ping-pong lemma.

\section{Paradoxical decomposition of group actions} 
Suppose that $G$ acts on $X$ and $\mathfrak{g}=(g_1,g_2,\ldots,g_n)$ is an ordered subset of $G$, and $\mathcal{E}=\{E_1,E_2,\ldots,E_m\}$ is a finite partition of $X$. By the concept of configurations for group action related to $(\mathfrak{g},\mathcal{E})$ we mean an $(n+1)$-tuple $C=(C_0,C_1,\ldots,C_n)$ such that  $C_i\in \{1,\ldots,m\}$, for each $i\in \{0,\ldots,n\}$ and there exists $x\in E_{C_0}$ with the property that $g_i \cdot x \in E_{C_i}$, for each $i \in \{1,\ldots,n\}$. We say that $(\mathfrak{g},\mathcal{E})$ is a configuration pair if $\mathfrak{g}$ is an ordered subset of $G$ and $\mathcal{E}$ is a finite partition of $X$.
We will denote the set of all configurations of group action related to $(\mathfrak{g},\mathcal{E})$ by $\mathrm{Con}(\mathfrak{g},\mathcal{E};X)$ (or $\mathrm{Con}(\mathfrak{g},\mathcal{E})$ again, if there is no ambiguity). Also we define 
\begin{equation}
\mathrm{Con}(G\curvearrowright X)=\{\mathrm{Con}(\mathfrak{g},\mathcal{E};X) : (\mathfrak{g},\mathcal{E})\;\text{is a configuration pair}\}.
\end{equation}
If $G=\langle\mathfrak{g}\rangle=\langle g_1,\ldots,g_n \rangle$ acts on $X$ and $\mathcal{E}=\{E_1,\ldots,E_m\}$ is a partion for $X$, then 
$C=(C_0,\ldots,C_n)\in \mathrm{Con}(\mathfrak{g},\mathcal{E};X)$ if and only if $x_0(C)\neq \emptyset $, where   $x_0(C)=E_{C_0}\cap g_1^{-1}E_{C_1}\cap\ldots\cap g_n^{-1}E_{C_n}$.\\
Let $G\curvearrowright X$ be a group action and $(\frak g, \mathcal E)$ be as above. Then it is not difficult to see that $\{x_0(C): \ C\in Con(\frak g, \mathcal E)\}$ forms a partition for $X$ (The best general reference here is \cite{rw}). Since for each $j$, $g_j\cdot X=X,$ it follows that the families $\{x_j(C): \ C\in Con(\frak g, \mathcal E)\}$  make partitions for $X$, as well. It is also easy to see that 
\begin{equation}\label{efraz} 
E_i=\bigsqcup_{C_j=i}x_j(C),\qquad 1\leq i\leq m.
\end{equation}

Designate $f_C$ to every configuration $C$, then
\begin{equation}
\sum \{f_C : C_j=i\}=\sum \{f_C : C_0=i\}\qquad (C\in Con(\mathfrak{g},\mathcal{E};X))
\end{equation}
are called configuration equations related to configuration pair $(\mathfrak{g},\mathcal{E})$. We denote this system by $\mathrm{Eq}(\mathfrak{g},\mathcal{E};X)$. Any solution of $\mathrm{Eq}(\mathfrak{g},\mathcal{E};X)$ that satisfies $\sum \{f_C : C\in Con(\mathfrak{g},\mathcal{E};X) \}=1$ and $f_C\geq 0$ is called a normalized solution.

The statement of \cite[Proposition 3.2]{rw} remains true when we use configurations for group actions. We now give the proof with the appropriate adjustments in the notations. 
\begin{theorem}
Let $G$ be a finitely generated group acting on a set $X$. Then the following statements are equivalent 
\begin{itemize}
\item[(1)] $G\curvearrowright X$ is amenable,
\item[(2)] $\mathrm{Eq}(\mathfrak{g},\mathcal{E};X)$ has a normalized solution for each configuration pair $(\mathfrak{g},\mathcal{E})$.
\end{itemize}
\end{theorem}
\begin{proof}
We first suppose that $G$ has an amenable action on $X$. Let $m$ be the left invariant mean on $l_{\infty}(X)$ and $(\mathfrak{g},\mathcal{E})$ a configuration pair. Let $y_C:=\langle \chi_{x_0(C)},m\rangle$ for the configuration $C\in\mathrm{Con}(\mathfrak{g},\mathcal{E};X)$. Then $y_C\geq 0$ and 
\begin{align*}
\sum\{y_C:\;C_j=i\} & =\left\langle\sum\left\{\chi_{x_0(C)}:\;C_j=i\right\},m\right\rangle\\
& =\langle\delta_{g_j}\cdot\chi_{E_i},m\rangle=\langle\chi_{E_i},m\rangle\\
& =\left\langle\sum\left\{\chi_{x_0(C)}:\;C_0=i\right\},m\right\rangle\\
&=\sum\{y_C:\;C_0=i\}.
\end{align*}
Also 
\begin{align*}
\sum\left\{\chi_{x_0(C)}:\;C\in\mathrm{Con}(\mathfrak{g},\mathcal{E};X)\right\} & =\sum_{i=1}^{m}\sum\left\{\chi_{x_0(C)}:\;C_0=i\right\}\\
& =\sum_{i=1}^{m}\chi_{E_i}\\
&=1.
\end{align*}
By the last equality
\begin{align*}
\sum\left\{y_C:\;C\in\mathrm{Con}(\mathfrak{g},\mathcal{E};X)\right\} & =\left\langle\sum\left\{\chi_{x_0}(C):\;C\in\mathrm{Con}(\mathrm{g},\mathcal{E};X)\right\},m\right\rangle\\
&=\langle 1,m\rangle\\
&=1.
\end{align*}
All in all it shows that $(y_C)_{C\in\mathrm{Con}(\mathfrak{g},\mathcal{E};X)}$ is a normalized solution for the configuration equations.\\
Conversely assume that $\{y_C:\;C\in\mathrm{Con}(\mathrm{g},\mathcal{E};X)\}$ is a normalized solution for $\mathrm{Eq}(\mathfrak{g},\mathcal{E};X)$. Choose $x_0^C\in x_0(C)$, where $C\in\mathrm{Con}(\mathfrak{g},\mathcal{E};X)$.\\
Define the function $f_{(\mathfrak{g},\mathcal{E})}$ on $X$ as follows
$$
f_{(\mathfrak{g},\mathcal{E})}(x):= \begin{cases}
y_C &  x=x_0^C\\
0 &  \text{otherwise}
\end{cases}.
$$
Normality of the solution guarantees that $f_{(\mathfrak{g},\mathcal{E})}\in P(X)=\{f\in l^1(X): f\geq 0\;, \parallel f \parallel_1 = 1\}$, $\left(f_{(\mathfrak{g},\mathcal{E})}\right)_C=y_C$ and by
\begin{equation*}
\langle f,\delta_{g_i}\cdot\chi_{E_i}-\chi_{E_i}\rangle=\sum\{f_C:\;C_j=i\}-\sum\{f_C:\;C_0=i\}
\end{equation*}
we can conclude that
\begin{equation*}
\langle f_{(\mathfrak{g},\mathcal{E})},\delta_{g_j}\cdot\chi_{E_i}-\chi_{E_i}\rangle=0\;.
\end{equation*}

This theorem in particular implies that if $G_1\curvearrowright X_1$ and $G_2\curvearrowright X_2$ are two group actions such that  $Con(G_1\curvearrowright X_1) =Con(G_2\curvearrowright X_2),$ the amenability of one of them yields the  amenability of the other one. 

Let $\mathcal E=\{E_1,\dots,E_m\}$ be a partition for a set $X$. Recall that a partition $\mathcal E'=\{E'_1,\dots,E'_t\}$ is a refinement for $\mathcal E$ if for each $i\in\{1,\dots,t\}$ there exists $l\in\{1,\dots,m\}$ such that $E'_i\subseteq E_l$. Suppose that $\mathfrak g=(g_1,\dots,g_n)$ is an ordered string of elements of a group $G$ and $G$ acts on $X$. It follows immediately that for any $C=(C_0,C_1,\dots,C_n)\in Con(\mathfrak g,\mathcal E')$ there exists a unique $D=(D_0,D_1,\dots,D_n)\in Con(\mathfrak g,\mathcal E)$ such that for each $j\in\{0,1,\dots,n\}$, $E'_{C_j}\subseteq E_{D_j}$.
In this case we use the notation $C\ll D$. 
All this facts show that  any $w^{\ast}$-limit point of $\left(f_{(\mathfrak{g},\mathcal{E})}\right)_{(\mathfrak{g},\mathcal{E})}$ is a left invariant mean on $l_{\infty}(X)$, and consequently the action of the group is amenable, which is our claim.
\end{proof}

%In \cite{arw} they showed that if $Con^*(G)\subseteq Con^*(\mathbb{F}_n)$ then $G\cong \mathbb{F}_n$. By a similar argument as is used in \cite{arw}, the following theorem is immidiate.
%\begin{theorem}
%Let $G\curvearrowright X$, $\mathcal{E}=\{E_i^{\pm}: 1\leq i \leq n\}\cup \{ E_1\}$ be a partition of $X$ and $G=\langle g_1,\ldots ,g_n \rangle$ such that the conditions 1 to 5 holds. Then $\mathbb{F}_n\leq G$ and conversely.
%\end{theorem}

Let $G=\mathbb{F}_2$ be the non-abelian free subgroup on two generators. In \cite{rw} they showed that $Eq(\mathfrak{g},\mathcal{E})$ for some generator $\mathfrak{g}$ of $G$ and a partition $\mathcal{E}$ of $G$, has no normalized solution. By a refinement $(\mathfrak{g}',\mathcal{E}')$ of $(\mathfrak{g},\mathcal{E})$, they constructed a paradoxical decomposition for $G$. In the following theorem we show that if $Eq(\mathfrak{g},\mathcal{E})$ has no normalized solution , then is the same for each refinement pair $(\mathfrak{g}',\mathcal{E}')$ of $(\mathfrak{g},\mathcal{E})$.
\begin{question}
Let $Eq(\mathfrak{g},\mathcal{E})$ has no normalized solution. Is there a refinement $(\mathfrak{g}',\mathcal{E}')$ of $(\mathfrak{g},\mathcal{E})$ such that $(A_i,B_j;g_i,h_j)$ is a paradoxical decomposition of $G$ where $\mathcal{E}'=\{A_i,B_j:1\leq i \leq n , 1\leq j \leq m \}$ and $\mathfrak{g}'=(g_1,\ldots,g_n;h_1,\ldots,h_m)?$ 
\end{question}
\begin{definition}
Let $(\mathfrak{g}',\mathcal{E}')$ and $(\mathfrak{g},\mathcal{E})$ be two configuration pairs for the group action $G\curvearrowright X$ such that $\mathcal E'$ is a refinement for $\mathcal E$, $\mathfrak g=(g_1,\dots,g_n)$ and $\frak g'=(g_1,\dots,g_n,g_{n+1},\dots,g_s)$.  Then $(\mathfrak{g}',\mathcal{E}')$ is called a refinement for $(\mathfrak{g},\mathcal{E})$.
\end{definition}
 In the following theorem we investigate the relation between the existence of a normalized solution for $Eq(\mathfrak g,\mathcal E;X)$ and $Eq(\mathfrak g',\mathcal E';X)$, when $(\mathfrak{g}',\mathcal{E}')$ is a refinement for $(\mathfrak{g},\mathcal{E})$ (This can be found in \cite{yyy} for the case of configuration of groups).

\begin{theorem} 
Let $G\curvearrowright X$ and $(\mathfrak{g}',\mathcal{E}')$ be a refinement for $(\mathfrak{g},\mathcal{E})$. If $Eq(\mathfrak g',\mathcal E';X)$ admits a normalized solution. Then so does $Eq(\mathfrak g,\mathcal E;X)$.
\end{theorem}
\begin{proof}
We proceed to the proof in three steps.

\textbf{Step 1.}  By definition of configurations of group action and the notation  $C\ll D$, it is proved that 
\begin{equation}\label{lll}
Con(\mathfrak g,\mathcal E';X)=\bigsqcup_{D\in Con(\mathfrak g,\mathcal E;X)}\left\{C\in Con(\mathfrak g,\mathcal E';X) :\qquad C\ll D\right\}.
\end{equation}
Let $(z_C)_{C\in Con(\mathfrak g,\mathcal E';X)}$ be a normalized solution for $Eq(\mathfrak g,\mathcal E',X)$. For $D\in Con(\mathfrak g,\mathcal E;X)$ set $z_D=\sum_{C\ll D}z_C$. Then by (\ref{lll})
\begin{align*}
\sum_{ x_0(D)\subseteq E_i}z_D &= \sum_{ x_0(D)\subseteq E_i}\sum_{C\ll D}z_C\\
&= \sum_{\substack{E'_p\in \mathcal E', \\ E'_p\subseteq E_i}} \sum_{ x_0(C)\subseteq E'_p}z_C \\ 
%{ E'_p\in \mathcal E', E'_p\subseteq E_i}\sum_{ x_0(C)\subseteq E'_p}z_C \\ 
&= \sum_{\substack{E'_p\in \mathcal E', \\ E'_p\subseteq E_i}} \sum_{ x_j(C)\subseteq E'_p}z_C\\&= \sum_{ x_j(D)\subseteq E_i}\sum_{C\ll D}z_C\\&=\sum_{ x_j(D)\subseteq E_i}z_D.
 \end{align*}
On the other hand one has
\begin{equation*}
\sum_{D\in Con(\frak g,\mathcal E;X)}z_D=\sum_{D\in Con(\frak g,\mathcal E;X)}\sum_{C\ll D}z_C=\sum_{C\in Con(\frak g,\mathcal E';X)}z_C=1.
\end{equation*}
Therefore $(z_D)_{D\in Con(\mathfrak g,\mathcal E;X)}$ is a normalized solution for $Eq(\mathfrak g,\mathcal E;X)$. 

\textbf{Step 2.}
Suppose that $Eq(\mathfrak g',\mathcal E;X)$ admits a normalized solution  $(z_C)_{C\in Con(\mathfrak g',\mathcal E;X)}$. Let $C=(C_0,C_1,\dots,C_n,C_{n+1},\dots, C_s)\in Con(\mathfrak g',\mathcal E;X)$. Then there exists an unique $D=(D_0,D_1,\dots,D_n)\in Con(\mathfrak g,\mathcal E;X)$  such that $$C_j=D_j,\qquad 0\leq j\leq n.$$ In this case we use the notation $D\preceq C$. For $D\in Con(\mathfrak g,\mathcal E;X)$  set 
$$z_D=\sum_{D\preceq C}z_C.$$
Then we have 
\begin{align*}
\sum_{\substack{ D\in Con(\frak{g},\mathcal E;X),\\ x_0(D)\subseteq E_i}}z_D &= \sum_{\substack{  D\in Con(\frak{g},\mathcal E;X),\\ x_0(D)\subseteq E_i}}\sum_{D\preceq C}z_C\\
&= \sum_{ \substack{ C\in Con(\frak{g'},\mathcal E;X),\\ x_0(C)\subseteq E_i}}z_C \\ 
&= \sum_{ \substack{ C\in Con(\frak{g'},\mathcal E;X), \\ x_j(C)\subseteq E_i}}z_C\\&= \sum_{\substack{ D\in Con(\frak{g},\mathcal E;X), \\ x_j(D)\subseteq E_i}}\sum_{D\preceq C}z_C\\&=\sum_{ \substack{ D\in Con(\frak{g},\mathcal E;X),\\ x_j(D)\subseteq E_i}}z_D.
 \end{align*}
Besides, we have
\begin{equation*}
\sum_{D\in Con(\frak g,\mathcal E;X)}z_D=\sum_{D\in Con(\frak g,\mathcal E;X)}\sum_{D\preceq C}z_C=\sum_{C\in Con(\frak g',\mathcal E;X)}z_C=1.
\end{equation*}
Thus  $(z_D)_{D\in Con(\mathfrak g,\mathcal E;X)}$ is a normalized solution for $Eq(\mathfrak g,\mathcal E;X)$. 

\textbf{Step 3.} If $Eq(\mathfrak{g}',\mathcal{E}';X)$ has a normalized solution, then by Step 1, $Eq(\mathfrak{g},\mathcal{E}';X)$ admits a normalized solution and by Step 2, so does $Eq(\mathfrak{g},\mathcal{E};X)$, and the proof is complete.
\end{proof}

The following two results are concerning the cardinal number of sets in configuration group actions.

\begin{lemma}
If $G\curvearrowright X$ and $H\curvearrowright Y$ are group actions and $\mathrm{Con}(G\curvearrowright X)=\mathrm{Con}(H\curvearrowright Y),$ then
\begin{itemize}
\item[(1)] $X$ is finite (resp. infinite) if and only if $Y$ is finite (resp. infinite).
\item[(2)] $|X|=n$  if and only if $|Y|=n$.
\item[(3)] $|X|\geq n$ if and only if $|Y|\geq n$. 
\end{itemize}
\end{lemma}
\begin{proof}
Suppose $|Y|> n=|X|$, and $Con(G\curvearrowright X)=Con(H\curvearrowright Y)$. Take $\mathcal{F}=\{F_1,\ldots F_{n+1}\}$, where $F_1 =\{y_1\},\ldots ,F_n=\{y_n\},F_{n+1}=Y\setminus (\cup _{i=1}^{n} F_i)$, for some $y_i\in Y$, and a generating set  $\mathfrak{h}$ for $H$. There exists a configuration pair $(\mathfrak{g} ,\mathcal{E})$ for $X$ such that $Con(\mathfrak{g},\mathcal{E})=Con(\mathfrak{h},\mathcal{F})$. Thus $n=|X|\geq |\mathcal{E}|=|\mathcal{F}|=n+1$, which is a contradiction.
%We begin with supposing that $X$ and $Y$ are finite and $|Y|\gslant |X|=n\in\mathbb{N}$. Let $y_1,\dots,y_n\in Y$ and define $E_1=\{y_1\}$, $E_2=\{y_2\}$, $\ldots$, $E_n=\{y_n\}$, $E_{n+1}=Y\setminus\cup_{i=1}^{n}E_i$.  Then $\mathcal{E}:=\{E_i:\;i\leq i\leq n+1\}$ makes a partition for $Y$. Assume that $\mathfrak{h}$ is an ordered string in $H$. Hence there exists a pair $(\mathfrak{g},\mathcal{F})$ such that $\mathfrak{g}$  is a string in $G$ and $\mathcal{F}$ is a partition of $X$ and besides, $\mathrm{Con}(\mathfrak{g},\mathcal{F};X)=\mathrm{Con}(\mathfrak{h},\mathcal{E};Y)$. Then $|\mathcal{E}|=|\mathcal{F}|\leq n$, which is a contradiction.
%Now suppose that $X$ is countable. Then for each $n\in\mathbb{N},\quad |X|>n$. This implies that $|Y|>n$ and so $Y$ is at least countable, and this is precisely the assertion of the lemma.

\end{proof}
\begin{corollary}
Let $G$ and $H$ be two finitely generated groups such that $G\curvearrowright X$ and $H\curvearrowright Y$. If these actions are transitive and $\mathrm{Con}(G,X)=\mathrm{Con}(H,Y),$ then $|X|=|Y|$. 
\end{corollary}
\begin{proof}
Let us first recall that the actions are transitive, so for any $x\in X$ and $y\in Y$, we have $|X|=|\frac{G}{Stab(x)}|$ and $|Y|=|\frac{G}{Stab(y)}|$. Then if $X$ and $Y$ are both infinite, the result is obvious because $G$ and $H$ are finitely generated and so countable. The rest of the statement follows from the pervious Lemma. 
\end{proof}

\begin{theorem}\label{th1}
Let $G\curvearrowright X$ and $H\curvearrowright Y$ be two group actions and $f$ be a function  from $X$ onto $Y$. Let $\phi:G\rightarrow H$ be an epimorphism  such that  
$$f(g.x)=\phi(g).f(x)\qquad (x\in X,\ g\in G).$$
Then $Con(H\curvearrowright Y)\subseteq Con(G\curvearrowright X).$
\end{theorem}
\begin{proof}
For any partition $\mathcal E=\{E_1,\dots,E_m\}$ of $Y$ set $\bar{\mathcal E}=\{{\bar E}_1,\dots,{\bar E}_m\}$, where $\bar{E_i}=f^{-1}(E_i)$. Then it is readily seen  that $\bar{\mathcal E}$ is a partition of $X$. For each ordered subset $\mathfrak g=(g_1,\dots,g_n)$ of $G$ we have $Con(\frak g, \mathcal E)=Con(\frak h,\bar{\mathcal E})$, where $\frak g=(g_1,\dots,g_n)$ and $g_j=\phi (h_j)$. Indeed, for $C=(C_0,\dots,C_n)\in Con(\frak g, \mathcal E)$ and $x_0\in x_0(C)$, there exists $y_0\in Y$ such that  $x_0=f(y_0)$ and we have 
$$f(h_jy_0)=\phi(h_j)f(y_0)=g_jx_0\in E_{C_j}\qquad 1\leq j\leq n.$$
It means that $ h_jy_0\in \bar{E_i}$. Thus $C\in Con(\frak h,\bar{\mathcal E})$. This completes the proof.
\end{proof}
The following corollary is obtained immediately.
\begin{corollary}
Let $G$ be a group acting transitively on a set $X$. Then $Con(G\curvearrowright X)\subseteq Con(G)$.
\end{corollary}
\begin{theorem}
Let $Con(H\curvearrowright Y)\subseteq Con(G\curvearrowright X).$ Then $\tau(G\curvearrowright X)\leq \tau(H\curvearrowright Y).$
\end{theorem}
\begin{proof}
Let $(\frak g, \mathcal E)$ be a configuration pair of the group action $G\curvearrowright X$ . Since $\{x_0(C)\ C\in Con(\frak g, \mathcal E)\}$ is a partition for $X$  and by (\ref{efraz}) for each $1\leq i\leq m,$ $E_i=\bigsqcup_{C_j=i}x_j(C)$.\\
%By a similar argument as is used in the proof of \cite[Theorem 3.8]{rty} for the case of group action.\\
If $(E_1,\ldots,E_n;E_{n+1},\ldots,E_{n+m};g_1,\ldots,g_n;g_{n+1},\ldots,g_{n+m})$ is a paradoxical decomposition for the action $G\curvearrowright X$ such that $Con(G\curvearrowright X)=Con(H\curvearrowright Y)$ and $\mathcal{E}=\{E_1,\ldots,E_{n+m}\}$, $\mathfrak{g}=\{g_1,\ldots,g_{n+m}\}$, then there exists a configuration pair $(\mathfrak{g}',\mathcal{E}')$ for the action $H\curvearrowright Y$ such that $Con(\mathfrak{g},\mathcal{E};X)=Con(\mathfrak{g}',\mathcal{E}';Y)$ and $(E'_1,\ldots,E'_n;E'_{n+1},\ldots,E'_{n+m};g'_1,\ldots,g'_n;g'_{n+1},\ldots,g'_{n+m})$ is a paradoxical decomposition for the action $H\curvearrowright Y$.
\end{proof}
\begin{corollary}
Two configuration equivalent group actions admit the same Tarski numbers. 
\end{corollary}

\begin{remark} 
There are several examples showing that in general $Con(G)\neq Con(G\curvearrowright X).$ In fact every group action for which $\tau(G)\neq \tau(G\curvearrowright X)$ is an example.
For a special case let $G$ be a group containing a non-abelian free group acting on a set $X$ and the point stabilizers of this action be all cyclic. Then $\tau(G)=\tau(G\curvearrowright X)=4$. Whereas if this group acts trivially on the same set $X$ (i.e. g.x=x), then $\tau(G\curvearrowright X)\neq 4,$ by \cite[Theorems 4.5 and 4.8]{wagon}. In this example, if $\frak g=(g_1,\dots,g_n)$ and $\mathcal E=\{E_1,\dots,E_m\}$ are arbitrary ordered string and partition of $G$, then
$$
Con(\frak{g}, \mathcal E;X)=\{(i,i,\dots,i),\qquad 1\leq i\leq m\}.
$$
So  the system of configuration equations $Eq(\frak{g}, \mathcal E;X)$ has normalized solution $z_{C}=1/m$, $C\in Con(\frak g,\mathcal E;X)$. Since $\frak{g}$ and $\mathcal E$ are chosen arbitrary, $\tau(G\curvearrowright X)=\infty$.
\end{remark}
\begin{remark}
Let the action $G\curvearrowright X$ be transitive and $x_0\in X$ and $Stab_{x_0}=\{g\in G, gx=x\}$ be the stabilizer of $x_0$ in $G$. Assume that $f:X\rightarrow \frac{G}{Stab(x_0)}$, is defined through $x\mapsto g_x Stab(x_0)$, where $g_x\in G$ is an element that $x=g_x.x_0$. Then $f$ is a $G$-invariant map onto $\frac{G}{Stab(x_0)}.$ Therefore $Con(G\curvearrowright X)\subseteq Con(G)$.
\end{remark}

\begin{corollary}
If $G\curvearrowright X$ is free transitive, then $Con(G\curvearrowright X)= Con(G)$.
\end{corollary}
%Let $G$ be a group and $N$ be a normal subgroup of $G$. $G$ acts on $\frac{G}{N}$ through the following action
%$$g.(xN)=gxN.$$ The following theorem, which is a direct consequence of Theorem \ref{th1}, gives the relation between the configuration sets of the action of $G$ on $G$ and $\frac{G}{N}$.
%\begin{theorem} 
%With the above notations $Con(G\curvearrowright \frac{G}{N})\subseteq Con(G).$
%\end{theorem}
\section{Configurations and ping-pong lemma}
In this section we state a number of the ping-pong lemma's statements. In a special case these statements ensure that several elements in a group acting on a set, freely generate a free subgroup of that group. Then we investigate the connection between configurations and  ping-pong lemma.

\begin{theorem}[Ping-Pong lemma]\cite{harp} \label{ping two}
Let $G$ be a group and let $X$ be a set equipped with a $G$-action. Suppose $\Gamma$ is generated by subgroups $\Gamma_1$ and $\Gamma_2$ with $|\Gamma_1| \geq 3$ and $|\Gamma_2| \geq 2$. Suppose furthermore that $X_1$ and $X_2$ are two subsets of $X$ with $X_2$ not included in $X_1$. Finally, suppose that for each nonidentity elements $\gamma_i\in \Gamma_i,\ \gamma_1(X_2)\subseteq X_1$ and $\gamma_2(X_1)\subseteq X_2$. Then $\Gamma\cong \Gamma_1 \ast\Gamma_2$.

\end{theorem}
The ping-pong lemma can be extended for more than two players. The proof of the next theorem is similar to the case of two subgroups (The classical work here is \cite[II. B ]{harp1}).

\begin{theorem}[Ping-Pong Lemma for several subgroups] \label{ping more}
Let $G\curvearrowright X$ be a group action and  $H_1,\ldots,H_k$ be subgroups of  $G$ such that $|H_i|>2$  for some $i\in\{1,\ldots,k\}$. Suppose that there are $k$ pairwise disjoint subsets  $X_1,\ldots,X_k$ of $X$ such that 
for each $i\neq s$ and for any  non-identity element $h\in H_i$, $h(X_s)\subseteq X_i$.
Then $\langle H_1,\ldots,H_k\rangle=H_1\ast\dots\ast H_k$.
\end{theorem}

For the case of cyclic subgroups, the ping-pong lemma has the following corollary. 

\begin{corollary}
Let $k\geq 2$ be an integer, $G\curvearrowright X$ and $A_1,\ldots,A_k,B_1,\ldots,B_k$  be $2k$ pairwise disjoint subsets of $X$ and  $g_1,\ldots,g_k\in G$ be such that 
\begin{equation*}
 B_i^c\subseteq g_iA_i\, \qquad 1\leq i\leq k.
\end{equation*}
Then $\langle g_1,\ldots,g_k\rangle$ is a free subgroup of $G$.
\end{corollary}
\begin{proof}
For $1\leq i\leq k$ put $X_i=A_i\cup B_i$. Then for $i\neq s$
$g_iX_s=g_i(A_s\cup B_s)\subseteq g_iA_i^c\subseteq B_i.$
and $g_iB_i\subseteq g_iA_i^c\subseteq B_i$. Therefor for each $n\in \mathbb N,$
$g_i^nX_s\subseteq X_i.$
On the other hand $$g_i^{-1}X_s=g_i^{-1}(A_s\cup B_s)\subseteq g_i^{-1}B_i^c\subseteq A_i$$ and $$g_i^{-1}A_i\subseteq g_i^{-1}B_i^c\subseteq A_i.$$ Hence for any positive integer number $n$, $g_i^{-n}X_s\subseteq X_i.$ Consequently for any non-identity element $h$ of $\langle g_i\rangle$, we have $hX_s\subseteq X_i$. Then by Theorem \ref{ping more} $\langle g_1,\dots,g_k\rangle\cong \langle g_1\rangle\ast\dots\ast\langle g_k\rangle.$ Clearly each $g_i$ is of infinite order. Therefore $\langle g_1,\dots,g_k\rangle\cong \mathbb F_k$, and the corollary follows.
\end{proof}

There are many papers in the literature, which use ping-pong lemma to prove the existence of a free subgroup. It is of real interest how this lemma relates to configuration sets. 

\begin{proposition} If $\mathrm{Con}(G)=\mathrm{Con}(H)$ and $G$ satisfies the conditions of  Lemma \ref{ping two}, then these conditions hold for $H$. 
\end{proposition}
\begin{proof}
By assumption $G$ contains the non-abelian free group $\mathbb F_2$. By proposition \ref{t4} $H$ contains $\mathbb F_2$ as well and this subgroup clearly admits the ping-pong lemma conditions.
\end{proof}

\begin{proposition}
The following statements are equivalent
\begin{itemize}
\item[(1)] $G$ is none-abelian;
\item[(2)] there exists $E_1,E_2,E_3,E_4,E_5$ of pairwise disjoint subsets of $G$ and $g_1,g_2 \in G$, which satisfies  $g_1E_1=E_2=g^{-1}_2 E_4$, and $E_3=g^{-1}_1 E_5=g_2 E_1$.
\end{itemize}
\end{proposition}
\begin{proof}
($(1)\Rightarrow(2)$) Since $G$ is not abelian, there exist $g_1,g_2\in G$ such that $g_1g_2\neq g_2g_1$. Actually by defining $E_1=\{e\}, E_2=\{g_1\}, E_3=\{g_2\}, E_4=\{g_1g_2\}, E_5=\{g_2g_1\}$, the result will be trivial.\\
($(2)\Rightarrow(1)$) Conversely, $g_1g_2E_1=g_1E_3=E_5$ and $g_2g_1E_1=g_2E_2=E_4$. Now since $E_4\cap E_5=\emptyset$, we can conclude that  $g_1g_2E_1\bigcap g_2g_1E_1=\emptyset$ and immediately $g_1g_2\neq g_2g_1$, thus $G$ is not abelian.
\end{proof}
\begin{lemma}
The following statements are equivalent
\begin{itemize}
\item[(1)] $a\in G$ is of infinite order.
\item[(2)] There exist disjoint sets $E_1$, $E_2$ such that $aE_1\subseteq E_1$, $aE_2\cap E_1\neq\emptyset$.
\end{itemize}
\end{lemma}
\begin{proof}

($(1)\Rightarrow(2)$) Let $E_1:=\{a^n:\;n\geq0\}$ then $aE_1\subseteq E_1$ is trivial, since $a^{-1}\in E_2$ thus $e\in aE_2\cap E_1$ and consequently $aE_2\cap E_1\neq\emptyset$.\\
($(2)\Rightarrow(1)$) On the contrary, suppose that $x\in aE_2\cap E_1$, yhus there exists $y\in E_2$ such that $x=ay$. Let $n\in \mathbb{N}$ be arbitrary, then $a^n y=a^{n-1} x \in E_1$. Now $a^n\neq e$ since otherwise $y\in E_1$, which is in contradiction with $E_1\bigcap E_2 = \emptyset$.
\end{proof}

We end this section with a version of ping-pong lemma to construct a paradoxical decomposition for a group or group action.

\begin{theorem}
Let group $G$ act on a set $X$ and there are subsets $X_1,\dots,X_n\subseteq X$ and $h_1,\dots,h_n\in G$ such that  $h_iX_i\subseteq X_{i+1},$ $X_{n+1}=X_1$ and $X=\bigcup_{i=1}^n (X_{i+1}\setminus h_iX_i).$ Then $G\curvearrowright X$ is paradoxical and $\tau(G\curvearrowright X)\leq n+2$.
\end{theorem}
\begin{proof}
Our proof starts with by setting $D_i=X_{i+1}\setminus h_iX_i$, $s_i=h_n h_{n-1}\dots h_i$, for $i\in \{1,\dots,n\}$ and $s_{n+1}=e$. Then by an inductive process one can see 
$$X_1=s_1X_1\sqcup s_2D_1\sqcup s_3D_2\dots\sqcup s_{n}D_{n-1}\sqcup s_{n+1}D_n.$$
We are now ready to construct a paradoxical decomposition for the action. Let $E_0=s_1X_1$ and $E_i=s_iD_i$, for $1\leq i\leq n$.  Recall that $X=X_1\sqcup {X_1}^c$, where ${X_1}^c$ is the complement of $X_1$ in $X$; i.e. $X=\left(\bigsqcup_{i=0}^n E_i\right)\sqcup   {X_1}^c$. We have
$$X=s_1^{-1}E_0\sqcup X_1^c.$$
On the other hand   by assumption 
$$X=\bigsqcup _{i=1}^nD_i= \bigsqcup_{i=1}^ns_{i+1}^{-1}E_i.$$
We have constructed a paradoxical decomposition with $n+2$ pieces. Therefore $\tau(G\curvearrowright X)\leq n+2$, which completes the proof.
\end{proof}

Recall from \cite{wagon}  that if $G$ is a group containing the non-abelian free group $\mathbb F_2$ and $G$ acts on a set $X$, then $\tau(G\curvearrowright X)=4$ if and only if the point stabilizers of this action are all cyclic. 

\begin{theorem}
Let group $G$ act on a set $X$ and there are subsets $X_1,X_2\subseteq X$ and $h_1,h_2\in G$ such that  $h_1X_1\subseteq X_{2}$ and $h_2X_2\subseteq X_{1}$  and  and $X=(X_{1}\setminus h_2X_2)\cup (X_{2}\setminus h_1X_1).$ Then $\tau(G)=4$ and the point stabilizers of this action are all cyclic.
\end{theorem}
\begin{proof}
We see at once that this theorem is the special case of the pervious theorem for $n=2.$ It means that $\tau(G\curvearrowright X)\leq 4$ and therefore $\tau(G\curvearrowright X)=4.$ Using \cite[Theorems 4.5 and 4.8]{wagon} $\mathbb F_2\subseteq G$ and the stabilizers of elements of $X$ are cyclic.
\end{proof}

\end{document}